\documentclass[amsart,11pt]{article}

\usepackage{amsfonts}
\usepackage{amsmath}
\usepackage{amssymb}
\usepackage{amsthm}
\usepackage{bbm}
\usepackage{graphicx}
\usepackage{gensymb}
\usepackage{comment}
\usepackage{placeins}
\usepackage{mathtools}
\usepackage{hyperref}
\usepackage{color}
\usepackage{algorithm} 
\usepackage{algorithmic}
\usepackage{authblk}
\usepackage[font=small,skip=0pt]{caption} 

\usepackage{lipsum}
\let\OLDthebibliography\thebibliography
\renewcommand\thebibliography[1]{
  \OLDthebibliography{#1}
  \setlength{\parskip}{0pt}
  \setlength{\itemsep}{0pt plus 0.3ex}
}

\usepackage{geometry}
\newcommand{\Z}{\mathbb Z}

\newcommand{\C}{\mathbb C}
\newcommand{\R}{\mathbb R}

\author[1]{Toby Sanders}
\affil[1]{School of Mathematical and Statistical Sciences, Arizona State University, Tempe, AZ, USA.}
\date{}
\title{Phase Based Alignment and Improved Projection Matching of Parallel Beam Tomography Data}
\begin{document}
\maketitle

\begin{abstract}
Tomography is an imaging technique that works by reconstructing a scene from acquired data in the form of line integrals of the imaging domain.  A fundamental underlying assumption in the reconstruction procedure is the precise alignment of the data values, i.e. the relationship between the data values and the paths of the lines of integration is accurately known.  In many applications, e.g. electron and X-ray tomography, it is necessary to establish this relationship using software alignment techniques or image registration due to misalignment when rotating the physical specimen.  Unfortunately, highly accurate software alignment is still a challenge to achieve in many cases, and improper alignment results in severe loss in the imaging resolution.  In this article, we develop a new approach that considers the alignment problem through a completely different lens, as a problem of recovering phase shifts in Fourier domain {within} the reconstruction algorithm.  The recovery of these phase shifts serves as the data alignment, which is done by calculating discrepancies between the misaligned data and the current reconstruction.  In the development of the approach, we investigate proper selection of parameters, and we show that it is fairly flexible and surprisingly accurate.  Finally, the analysis of our approach provides insight into why projection matching alignment by cross-correlation can be improved through low pass filtering, which we demonstrate.  Our methods are validated in a wide range of examples and settings.    

\end{abstract}

\section{Introduction}
Tomography is a tool for imaging two and three dimensional structures at a large range of scales, e.g. medical CT imaging, nanoscale electron tomography, and subsurface seismic tomography to name only a few.  Generally speaking, tomography is a rigorously established mathematical and scientific technique \cite{Natterer,Frank2006,kak2001principles,baruchel2000x,nolet1987seismic,frank2008electron}.  However unique challenges are inherent with each specific tomographic application, where the theoretical underpinnings of the mathematics do not align perfectly with the data acquisition mechanisms.  Therefore research is ongoing in tomography as such discrepancies are properly identified and new technologies push the boundaries of tomography into new applications \cite{hayashida2016practical,bouma2009fourier,midgley2009electron}.

This article addresses the preprocessing problem of data registration, which is inherent to a number of tomographic applications.  This includes electron tomography \cite{brandt2001automatic,amat2008markov,winkler2006accurate}, soft X-ray tomography \cite{parkinson2012automatic,mayo2007software}, synthetic aperture radar\footnote{Formally synthetic aperture radar is not a tomography problem, but can be reformulated as complex valued tomography problem using the Fourier slice theorem.}\cite{sanders2017combination,wahl1994phase}, and perhaps the most severe case occurring single particle cyro-electron microscopy \cite{frank2009single,penczek1994ribosome}. Loosely speaking, the imaging data comes in the form of line or local integrals of a function of the imaging object (e.g. the density), which are then processed to reconstruct an approximation of this function.  The data registration problem is that the location or path of each local integral from which the data is modeled may not be known precisely due to the nature of the data acquisition procedure.  In this setting, the data is usually preprocessed using data or image registration techniques before it is used to reconstruct the function, and even then it may be unknown whether the data registration procedure yielded the proper alignment.  Moving forward, in various applications improved sensing technologies may have the potential to alleviate this issue.  Nevertheless the issue is likely to persist as other imaging applications arise and the imaging limits are challenged, and image processing techniques will continue to be a necessary and powerful tool to overcome these problems.

The most common registration strategy is to simply cross correlate the neighboring data points and in hopes that globally the data all align properly \cite{frank1992alignment,castano2007fiducial,leary2013compressed,gursoy2017rapid}.  Unfortunately the cross-correlation model is purely heuristic.  The most obvious error resulting from cross correlation is cumulative drift in the registration, which means that at nearby perspectives the relative locations of the data values are well aligned, but small errors accumulate resulting in poor global alignment \cite{saxton1984three,brandt2001automatic}.  Although not the point of this article, we show several examples where cross correlation of a perfectly aligned data set results in a very poor reconstruction.

Another strategy is to align the data based on identification of bright features within the scene \cite{kremer1996computer,masich2006procedure,brandt2001automatic}.  These methods can generate very accurate results, particularly if there are many such features well distributed throughout the scene.  Unfortunately such features are not always available.  There is also the method of common lines, which attempts to align the projections over their common lines in the Fourier domain \cite{liu1995marker}.  More recently, center-of-mass based approaches have been developed \cite{sanders2015physically}, which essentially makes use of the \emph{Helgason-Ludwig
consistency conditions} \cite{helgason}.  To our knowledge, these methods are the first in which the registration relies on mathematical models that can be proven to yield exact results under certain conditions.  The conditions under which exact registration is theoretically achieved relies on first reducing the alignment to only translational shifts, and identifying certain cross-sections of the imaging scene in which the total mass of the projected scene is conserved throughout the full data sequence.  This cannot be guaranteed with every data set, but some additional developments have been made to improve further on these methods \cite{sanders2017mm}.

In this article we propose a new approach, which we refer to as phased based alignment (PBA).  First, the problem of unregistered data is reformulated into a problem of multiplicative phase angle errors in the data using a Fourier transform.  With this reformulation, we consider correcting for the data registration errors by equivalently estimating the phase shifts in the Fourier domain.  We estimate these phase shifts by calculating differences between the data and projections of the reconstructed image and iteratively update the registration in conjunction with the iterative reconstruction algorithm.  To assist in recovering the correct phase shifts, we make use of regularization methods for the reconstruction (e.g. nonnegativety or total variation), since for underdetermined problems methods such as filtered backprojection or ordinary least squares may recover a solution exactly matching the misaligned data values (hence little refinement may occur).  

Finally we mention an iterative alignment method that goes by several names, which we will refer to as projection matching \cite{dengler1989multi,penczek1994ribosome,parkinson2012automatic}.  We show this approach actually has some relation to the method proposed in this article.  We will address this later on in section \ref{sec:pm} and discuss why our method is an improvement to the basic projection matching.  Moreover, motivated by the analysis of our new approach, we demonstrate how one may notably improve the basic projection matching and achieve results comparable to our phase based approach.

The remainder of this article is organized in the following way.  We define the problem and motivate our method in section 2.  The exact details of the procedure along with some initial results and sample studies are given throughout section 3.  Finally, a thorough set of simulations are given in section 4, and in section 5 the method is applied to an electron tomography data set.  

\section{Problem Description and Reformulation}
We will characterize the problem for reconstruction of 2D cross-sectional images, which we will naturally extended to 3D later on.  Let $f(x,y)$ denote the image we want to reconstruct over some finite domain $\Omega$.  Then the ideal projection data take the form of the Radon transform
\begin{equation}\label{radon}
 Rf(x,\theta) = \int\displaylimits_{y: (x,y)Q_\theta^T\in \Omega} f((x,y)Q_\theta^T) \, dy,
 \end{equation}
  $$
 Q_\theta = \left[ \begin{array}{cc}
             \cos \theta & -\sin \theta \\
             \sin \theta & \cos \theta
            \end{array}
	    \right].
$$
In practice, $\theta$ will take on a discrete set of values, say $\{\theta_m\}_{m=1}^M$.  The problem of misalignment generally means that the data is not precisely acquired according to the model in (\ref{radon}), so that of the relative position of $f$ varies in an unknown way during data acquisition.  While this problem could manifest in several ways, for 2D we consider misaligned data of the form
\begin{equation}\label{radon-mis}
 \widetilde{Rf}(x,\theta_m) = Rf(x-\epsilon_m , \theta_m), \quad \text{for} \quad m=1,2,\dots , M.
\end{equation}


Proceeding forward, it will be necessary to characterize the problem in the discrete setting.  To do this we fix an $N$ point discretization over the integers and define the vectors $\tilde r_{\theta_m} = \{\widetilde{Rf}(n , \theta_m ) \}_{n=1}^N$ and $ r_{\theta_m} = \{{Rf}(n , \theta_m) \}_{n=1}^N$.  At this moment it is convenient to define a 2D sinogram matrix as
\begin{equation}\label{sinogram-def}
S\left( f,\{\theta_m \}_{m=1}^M \right)  = \left[ r_{\theta_1}  , r_{\theta_2}, \dots , r_{\theta_M} \right] ,
\end{equation}
i.e. each column of $S$ holds the Radon data at a particular angle.  Where implied, the sinogram may alternatively contain the misaligned data in the columns.

In this article we develop a new strategy for misalignment correction with a procedure that aligns based on the reconstruction.  Consider first if we wanted to reconstruct $f$ with the data while simultaneously determining \emph{optimal} registration shifts, e.g. determining a permutation matrix that aligns the data.  Taking Fourier transforms the misalignment error becomes a linear phase shift error.  

To describe this in more detail, let $g$ be some vector in $\C^N$.  Then the unitary discrete Fourier transform (DFT) of $g$ at frequency $k$ is given by
\begin{equation}\label{DFT}
 F(g)_k = \frac{1}{\sqrt{N}}\sum_{n=1}^N g_n e^{-i\frac{2\pi}{N}(k-1)(n-1)},
\end{equation}
where $k\in[1,N+1)$ and is only a non integer in the case of oversampling\footnote{We note that many other conventions could be used here.}.  If $g$ is offset like the projections, so that $\tilde g_n = g_{n+\epsilon}$, then
$
 F(\tilde g)_k = e^{i\frac{2 \pi}{N}\epsilon(k-1)} F(g)_k.
$
This equation represents a well known relationship where a shift of $g$ in the original domain results in a linear phase shift in the DFT of $g$.  Likewise, using the vectors $\tilde r_{\theta_m}$ and $ r_{\theta_m}$ as defined above and taking the corresponding DFT's yields
\begin{equation}\label{DFT-data}
 F(\tilde r_{\theta_m} )_k = e^{i\frac{2 \pi}{N} \epsilon_m (k-1)} F( r_{\theta_m} )_k.
\end{equation}
Hence the problem of recovering alignment shifts may be reformulated as recovering phase shifts in the Fourier domain.

Before arriving at our approach, let us first develop a few additional relationships.  First suppose that the phase shift $\frac{2 \pi}{N} \epsilon_m (k-1)$ in (\ref{DFT-data}) is in the interval $[-\pi, \pi)$, or equivalently that the shift is such that $|\epsilon_m | \le \frac{N}{2(k-1)}$.  Then using the complex logarithm (or alternatively the arc$\tan$) we find the shift given by
\begin{equation}\label{log2}
\epsilon_m =  \frac{N}{i2\pi (k-1)}\log\left( \frac{F(\tilde r_{\theta_m} )_k }{F( r_{\theta_m} )_k} \right) .
\end{equation}
On the other hand, for larger values of $k$ and $\epsilon_m$, we may have $\frac{2 \pi}{N} \epsilon_m (k-1) \notin [-\pi,\pi)$.  In this case, performing the same calculations as above we receive
\begin{equation}\label{log4}
 \tilde \epsilon_m = \epsilon_m + \frac{N}{k-1} \alpha =  \frac{N}{i2\pi (k-1)}\log\left( \frac{F(\tilde r_{\theta_m} )_k }{F( r_{\theta_m} )_k} \right).
\end{equation}
where $\alpha$ is the unique integer satisfying $\left( \frac{2 \pi}{N} \epsilon_m (k-1) + 2\pi \alpha \right) \in [-\pi ,\pi ).$ 

Finally, consider the effect of noise in these calculations.  We then redefine the acquired projections to include noise as 
\begin{equation}\label{radon-mis-2}
 \widetilde{Rf}(x,\theta_m) = Rf(x-\epsilon_m , \theta_m) + \eta_m(x) ,
\end{equation}
where $\eta_m$ is the noise term, and likewise for the discretized versions $\tilde r_{\theta_m}$ and $r_{\theta_m}$.  Suppose that the vectorized $\eta_m$ contains independent entries that are normally distributed with mean 0 and variance $\sigma^2$.  Then it can be shown that the unitary DFT of $\eta_m$ is also mean zero Gaussian distributed with variance $\sigma^2$ in the complex plane.  Then performing the series of operations as before, we get the estimated shift value as
\begin{equation}
 \tilde \epsilon_m =  \frac{N}{i2\pi (k-1)}\log\left( e^{i\frac{2\pi}{N} \epsilon_m (k-1)} + \frac{F(\eta_m)_k}{F(r_{\theta_m})_k} \right) .
\end{equation}

If $k$ is chosen to be small, then $F(r_{\theta_m})_k$ can be assumed to be significantly larger than $\sigma$ so that $| F(\eta_m)_k| \ll | F(r_{\theta_m})_k | $.  Under these safe assumptions, we see that the first term in the logarithm dominates, and we receive an accurate estimate of the shift.  Moreover, we will have many small values of $k$ to make use of that will allow us to \emph{average} out this noise.

We summarize these relationships with the following remarks:
\begin{enumerate}
 \item Obviously in practice we cannot determine $\epsilon_m$ by using (\ref{log2}) since  the true data $r_{\theta_m}$ used in the denominator is essentially what we're looking for to improve the reconstruction.  
 However, intermediate projections coming from approximate reconstructions from the misaligned data may allow us to estimate these values.
 
 \item For larger frequency values $k$, noise will play a bigger role and we may also receive wrapped phases.  Thus smaller values of $k$ will be more informative.
\end{enumerate}

\begin{figure}
\centering
 \includegraphics[width=.5\textwidth]{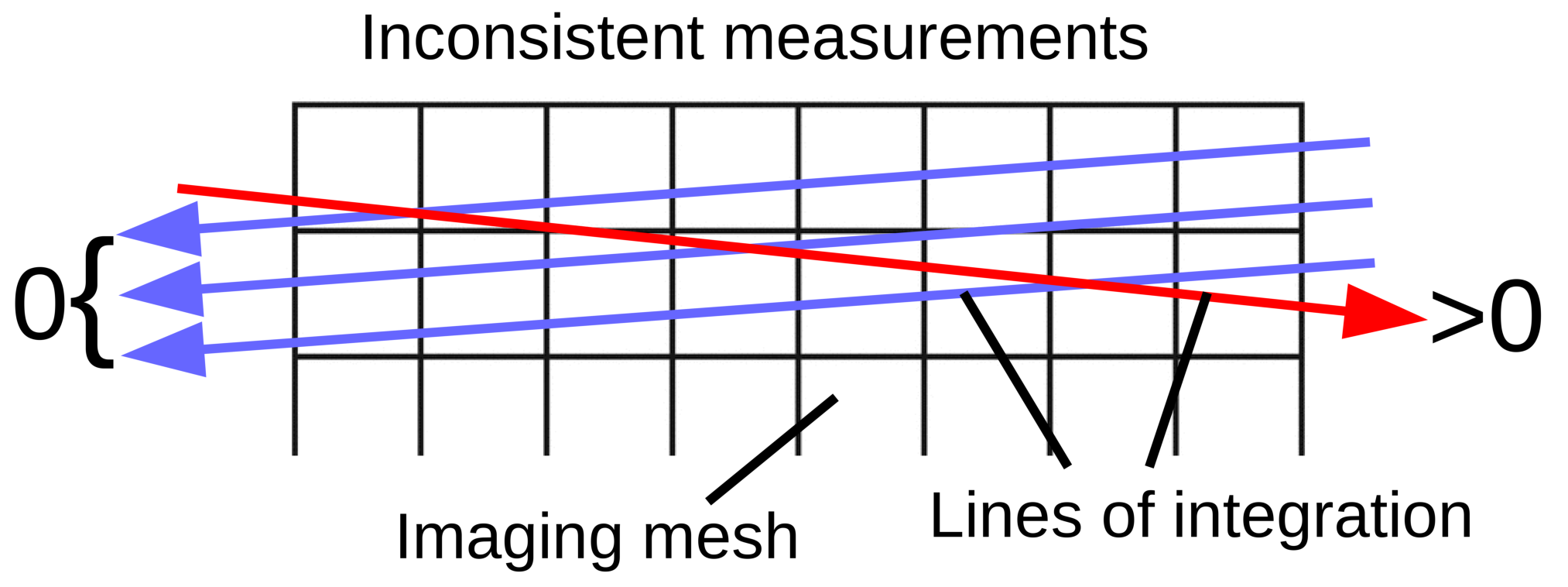}
 \caption{An illustration of contradictory measurements that can inform our phase based alignment.  If $f\ge0$, then having the blue lines of integration yield 0 and the red line yield a value greater than 0 is a contradiction.}
 \label{fig1}
\end{figure}

\section{Phase Based Alignment (PBA)} \label{sec:AF}
Let us use $f^{(j)}$ to denote some current reconstruction at iterate $j$ of $f$ from $\{ \widetilde r_{\theta_m} \}_{m=1}^M$.  Assuming the misalignment is not overwhelmingly large, we can infer that the projections of $f^{(j)}$ will satisfy to some extent the unknown ideal data, $ \{ r_{\theta_m} \}_{m=1}^M$.  This is especially true if we implement some regularizations for reconstructing $f^{(j)}$, such as nonnegativity or $\ell_1$, since misaligned data will result in inconsistent or contradictory data values realized through the regularization.  For example, recall that each data point is modeled as an integral of $f$ along some line.  Therefore, if some data point is zero, under a nonnegativity condition we deduce that $f$ is obviously zero along all points on the line.  However, misaligned data may give us inconsistencies in this regard (see Figure \ref{fig1}).  The misaligned data combined with the regularizations may result in inconsistencies throughout the entire image, resulting in the projections of $f^{(k)}$ satisfying to some extent both $\{ \widetilde r_{\theta_m} \}_{m=1}^M$ and $\{ r_{\theta_m} \}_{m=1}^M$.  On the other hand, methods such as ordinary least squares with no regularization may recover a solution which agrees exactly with the data, particularly for underdetermined problems, therefore motivating the use of regularization such as nonnegativity and total variation.

To this end, let us denote our forward projected solution of $f^{(j)}$ into the data space by $\{  r_{\theta_m} (f^{(j)}) \}_{m=1}^M$.  Then as computed in (\ref{log2}), we can now define an estimate for shifts by 
\begin{equation}\label{main-estimate}
\tilde \epsilon_m (k,j) =  \mathbf{Re}\Big\{ \frac{N}{i2\pi (k-1)}\log\left( \frac{F(\tilde r_{\theta_m} )_k }{F( r_{\theta_m} \left( f^{(j)} \right) )_k} \right) \Big\} ,
\end{equation}
where we take the real part since some complex parts of the expression may exist if the magnitudes of $F(\tilde r_{\theta_m} )_k$ and $F( r_{\theta_m} \left( f^{(j)} \right) )_k$ are not the same\footnote{We could also not divide by $i$ in (\ref{main-estimate}) and instead take the negative of the imaginary part.}.  In addition, noise in the data will inevitably be included in the calculations.
Notice the shift estimate depends on two parameters: the frequency $k$ used to estimate the shifts and the iteration(s) $j$ at which we choose to update the shifts.  The suitable choice of these parameters from which we then determine some \emph{ideal} shifts will be an important consideration.  For now, let's suppose from (\ref{main-estimate}) we have determined some ideal shift corrections $\{ \epsilon_m^* \}_{m=1}^M$.  Then we simply apply these shifts to the $\{ \tilde r_{\theta_m} \}_{m=1}^M$ and denote the updated data by $\{ \tilde r_{\theta_m}^{(j)} \}_{m=1}^M$.  We may then proceed with the reconstruction of $f^{(j+1)}$ from the updated data.  We now present the formal algorithmic description and some initial results.


\begin{figure*}
\centering
 \includegraphics[width=1\textwidth]{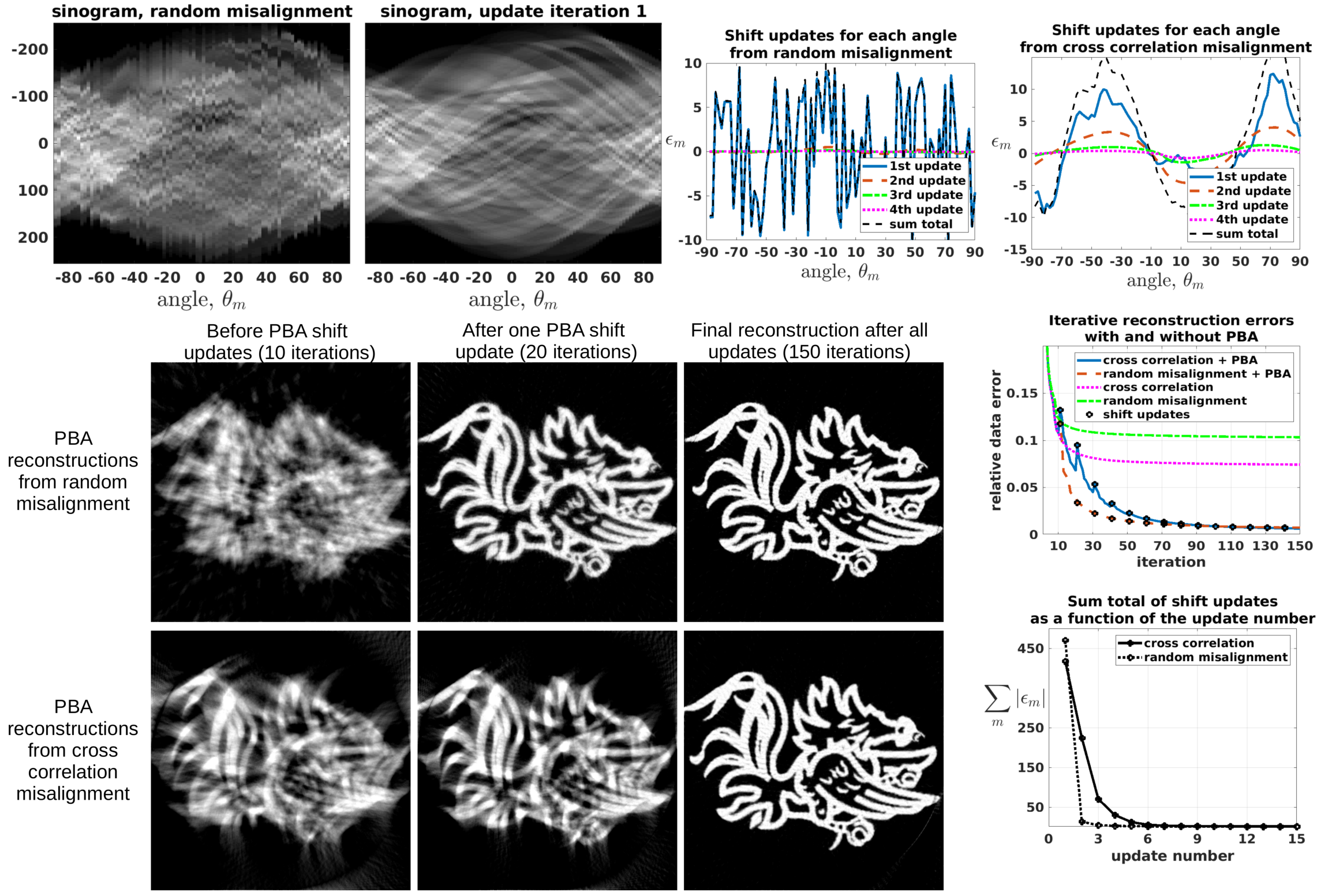}
 \caption{Phase based alignment results from cross correlation and random misalignment.  Top left: sinograms from the random misalignment and the PBA sinogram after one update.  Bottom left: iterative results of the PBA as the sinogram is updated for both random misalignment and misalignment due to cross correlation.  Top right: recovered phase based alignment shifts for each update for both random misalignment and cross correlation.  Bottom right: iterative reconstruction errors for with and without PBA (top) and sum total of the shift updates as a function of the update number (bottom). }
 \label{fig3}
\end{figure*}

\begin{figure}[ht]
\centering
 \includegraphics[width=.5\textwidth]{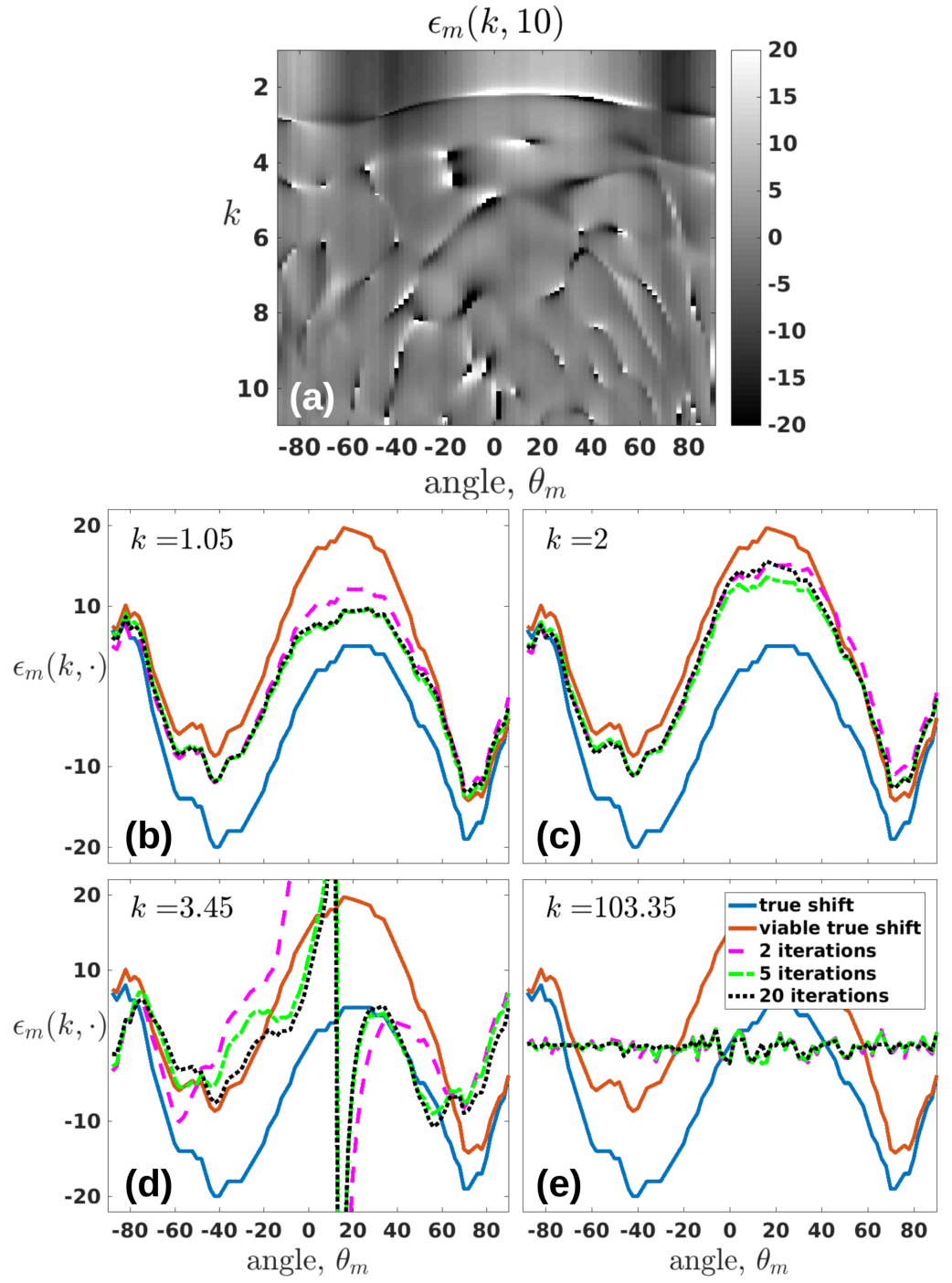}
 \caption{Demonstration of sinogram misalignment due to cross correlation, and determining shift values $\epsilon_m$ by (\ref{main-estimate}) for different $k$ values.  (b-e) Values of $\epsilon_m$ from (\ref{main-estimate}) for various $k$ and at different iterations.  (a) Values of $\epsilon(k,10)$ for many values of $k$.}
 \label{fig2}
\end{figure}

\subsection{Algorithm}

The general problem we've devised finds both $f$ and phase corrections in an alternating fashion.  Let $F$ denote the 1D discrete Fourier transform, and let $A_{\theta_m}$ denote the discrete Radon transform that maps $f$ to $r_{\theta_m}$.  Then let $ E(\epsilon_{m}) $ denote the diagonal matrix containing the phase corrections $\epsilon_{m}$ given by

\begin{equation}
 E(\epsilon_{m}) = 
 \left( \begin{array}{ccccc}
         e^{i\frac{2\pi}{N}(0)\epsilon_{m}} & 0  & \dots & 0\\
         0 &  e^{i\frac{2\pi}{N}(1)\epsilon_{m}}  & \dots & 0\\
         \vdots & \vdots  & \ddots & \vdots \\
         0 & 0 & \dots &  e^{i\frac{2\pi}{N}(N-1)\epsilon_{m}} 
 \end{array}\right) ,
\end{equation}
which will act as the data registration factor.

Then the objective minimization problem can be written as
\begin{equation}\label{eq:min1}
 \min_{f, \{\epsilon_{m}\} \in \R } \Big\{ \sum_{m=1}^M \| A_{\theta_m} f - 
 F^{-1}  E(\epsilon_{m}) F \tilde r_{\theta_m} \|_2^2 \Big\}
\end{equation}
Implementing the phase corrections by transforming into the Fourier domain by $F^{-1}  E(\epsilon_{m}) F \tilde r_{\theta_m}$, essentially just results in a circular shifting of $\tilde r_{\theta_m}$ by $\epsilon_{m}$.  In other words, this action acts as the data registration. These corrections are determined using (\ref{main-estimate}) for low values of $k$ and at some particular iterations $j$, which we determined at every 10th iteration to be suitable (see next section for justification), although this is flexible.  We can easily add in additional regularization terms suggested earlier to form a minimization such as 
\begin{equation}
\label{eq:min2}
 \min_{f, \{\epsilon_{m}\} \in \Z } \Big\{ \sum_{m=1}^M \| A_{\theta_m} f - 
 F^{-1}  E(\epsilon_{m}) F \tilde r_{\theta_m} \|_2^2 + H(f) \Big\}  \quad \text{s.t.} \quad f\ge0 ,
\end{equation}
where $H(f)$ is an image prior such as the TV norm.  The TV minimization in this paper uses the isotropic variant of the TV norm, which in its continuous formulation may be written as $\int_{\Omega} |\nabla f(x,y) | \, dx \, dy$ (for other variants and discretization, see for example \cite{bregman}).  The TV minimization is computed using the alternating direction method of multipliers (ADMM), which introduces a splitting variable and minimizes an augmented Lagrangian penalty function.  A pseudo algorithm is given in Algorithm \ref{algo-AF}.
\begin{algorithm}[!ht]
\caption{}
\label{algo-AF}
\begin{algorithmic}[1]
\STATE{Input data $\{ \tilde r_{\theta_m } \}_{m=1}^M$, and set iteration counts $J, L$:}
\STATE{Initialize $f^{(0,J)}$ and $\{ \tilde r_{\theta_m}^{(1)}\} = \{ \tilde r_{\theta_m } \}$.}
\FOR{$\ell=1$ \TO L}
\STATE{Set $f^{(\ell , 0)} = f^{(\ell-1,J)}$.}
\FOR{$j=1$ \TO J}
\STATE{Set $f^{(\ell,j)} = \Psi\left(f^{(\ell,j-1)}, \{\tilde r_{\theta_m}^{(\ell)}\} \right)$ , where $\Psi$ is some iterative updating procedure such as a gradient decent method or TV minimization.}
\ENDFOR
\STATE{Calculate shifts estimates $\tilde \epsilon_m(k,\ell) = \mathbf{Re}\Big\{ \frac{N}{i2\pi(k-1)} \log \left(\frac{F(\tilde r_{\theta_m} )_k}{F(r_{\theta_m} (f^{(\ell,J)}) )_k } \right) \Big\} $, for chosen $k$ and for $m=1,\dots,M.$}
\STATE{For each $m$, determine shifts $\epsilon_m(\ell) = \Phi(\{ \tilde \epsilon_m(k,\ell)\}_k)$, where $\Phi$ is some output or averaging function that uses (\ref{main-estimate}) and emphasizes the low values $k$.}
\STATE{Set $\{ \tilde r_{\theta_m}^{(\ell+1)} \}_m$ by circular shifting $\{ \tilde r_{\theta_m}^{(\ell)} \}_m$ by the values $\{ \epsilon_m (\ell)\}_m$.}
\ENDFOR
\STATE{Output final solution $f = f^{L,J}$ and final aligned data $\{r_{\theta_m}\}_m = \{\tilde r_{\theta_m}^{(L)}\}_m$}
\end{algorithmic}
\end{algorithm}

Demonstration of the algorithm is presented in Figure \ref{fig3}, where the results are presented for both random misalignment and misalignment brought on by cross-correlation of the perfectly aligned data\footnote{The original test image is given in Figure \ref{fig4}.}.  For the random misalignment, the misalignment shifts are independent integer values on the interval $[-10,10]$, each value with equal probability.  The suitable iteration count before updating the shift values (the value of $J$ in Algorithm \ref{algo-AF}) is set to 10.  The frequency values used for determining $\tilde \epsilon_m (k,j)$ are set to be $k\in (1,2]$ with an oversampling of 20, i.e. $k= 1.05, 1.10, \dots, 2.$  A detailed description and study for the determination of these parameters is given in section \ref{sec:param}, where we show the method is not particularly sensitive to the iteration counts between updates, but the chosen frequency values are very important.

Shown in the top left panel of Figure \ref{fig3}, using SIRT with a nonnegativity constraint, even after one phase based alignment (PBA) the sinogram is generally well aligned.  The reconstructions following this update are shown in the middle panel, label with 20 SIRT iterations.  This image still appears a bit \emph{fuzzy} or \emph{defocused}, hence some additional alignment may be necessary.  The final reconstruction after 14 sinogram alignment updates and 150 iterations of SIRT are shown to the right and indicate further improvement and near exact reconstruction.  In the bottom panel the same results are displayed for the cross correlation misalignment.  The bottom right panel shows the totals in the shift updates as a function of the update number.  That is $\sum_{m} |\epsilon_m|$ at the various update iterations.  As one would hope, it appears that the primary alignment corrections occur in the first few updates and converges after about 5 iterations.  The cross correlation requires more PBA iterations to converge, a result which has been previously observed \cite{mayo2007software,latham2016multi}.  In the middle right panel are the data errors (the normalized values of (\ref{eq:min1})) as a function of the number of SIRT iterations, and the results without PBA are also included.  Less surprising, the errors without the PBA corrections result in significantly higher values, hence inconsistencies in the data that cannot be sorted out in the reconstruction.  The total recovered shift values and iterative values are shown in the top right from both the random misalignment and cross correlation, and these also indicate that the PBA shifts converge after only a few iterations.

\subsection{Optimal Iteration Count and Frequencies}\label{sec:param}

As suggested in the early discussion and derivations in (\ref{log2})-(\ref{log4}), the smaller values of $k$ will be more helpful for determining the shifts, since avoiding phase wrapping in the logarithm in (\ref{main-estimate}) will be necessary.  It is not obvious to us how many iterations of the iterative reconstruction solver is appropriate.  For instance, we could argue one or two iterations is best in order to avoid over-fitting the misaligned data.  Alternatively, perhaps the algorithm should come close to convergence before we are able to realize accurate shift estimates.  Our numerical results presented in this section suggest that the results are not particularly sensitive to the number of iterations.  However, more iterations appear to give mildly more accurate shift estimates, but some sort of full convergence is certainly not necessary.  In general, if possible we would prefer get away with fewer iterations leading to faster convergence.  Our experience with this methods suggests that it is sufficient to update the shifts after every 5 to 10 iterations, and more importantly, we suggest the useful frequencies for (\ref{main-estimate}) should be $k\in (1,3)$.  

Figure \ref{fig2} highlights some of these points, where projection data is generated over a full $180\degree$ range on a 2D test image (the test image is shown in the top left of Figure \ref{fig4}).  The misalignment here is simply generated by using cross correlation from the perfectly aligned data.  The estimated shift values based on (\ref{main-estimate}) with the iteration count $j=10$ are visualized in (a) for increasing values of $k$.  Notice until around $k=3$ the shift estimates are essentially constant, at which point some values take a considerable jump.  This behavior persists with larger $k$ values, and these sorts of discontinuities over $k$ in the estimated shift values remain somewhat a mystery to us\footnote{Understanding these behaviors may provide useful information in the future.}.  In (b)-(e) these values along with the true shifts (and an alternative viable shift based off of the discussion in section \ref{sec:nonu})  are plotted for various values of $k$ and for various iteration counts $j=2,5,20$. Notice that the shift values are not particularly sensitive to the iteration count $j$, but again are sensitive to the frequency value $k$, suggesting for smaller values.

If not otherwise specified, for the purposed of this paper we choose to update the sinogram shifts after every 10th iteration, and we set the shifts by taking the average of $\epsilon_m(k,\cdot)$ for the values $k\in(1,2]$ with an over sampling of 20.  In other words we average the values for $k=1.05, 1.10, \dots, 2$.

\subsection{Extension to 3D Data Sets}\label{sec:3D}
The development of this work so far has been applicable to reconstruction of 2D objects from 1D projections or Radon data.  While for convenience research within tomographic imaging is often treated this way in mathematics, in many applications the goal is to reconstruct 3D objects from similar types of Radon data.  We outline this here, and show how the PBA can naturally be applied and possibly even improved in the 3D setting.  However, additional sources of misalignment may be present in 3D.

Let us continue denoting the structure or function by $f(x,y,z)$ defined over some domain $(x,y,z) \in \Omega \subset \R^3$.  We will define the Radon data by
\begin{equation}\label{3D-radon}
R f ( x,y, \theta) = \int\displaylimits_{z: (x,y,z) \in \Omega} f(x, (y,z) Q_\theta^T ) \, dz ,
\end{equation}
with $Q_\theta$ as in (\ref{radon}).  Hence, the data now are integrals of $f$ in the $z$ coordinate after rotation about the $x$-axis, and for a fixed $\theta$ we now have a whole 2D set of values over $x$ and $y$.  If we fix a cross-section say $x = x_0$, then it is clear to see that (\ref{3D-radon}) essentially reduces to the 2D definition in (\ref{radon}), where $(x,y)$ in (\ref{radon}) play the roles of $(y,z)$ in (\ref{3D-radon}).   As before we have $Rf$ for some set of angles $\{ \theta_m \}_{m=1}^M$.

The alignment problem in this case is that the data is given by 
\begin{equation}\label{3D-radon-unali}
\widetilde{Rf} (x,y, \theta_m ) = Rf (x - x_m , y -  \epsilon_m , \theta_m ) , 
\end{equation}
where there's now potentially an error in both the $x$ and $y$ directions.  We have purposefully labeled the error in the $y$ coordinate with $\epsilon_m$, since it will be analogous to the $\epsilon$ errors in the 2D setting.  We note here that there could be several other models and sources of misaligned data, such as offset tilt axis or inaccurate angle measurement.  Due to the nature of the method, we limit our work in this article to just translational shifts, and may address additional alignment issues in future work.  Offset tilt axis rotations can also be accounted for prior to our PBA \cite{wang2017fast}.

To reduce the 3D alignment problem to a similar problem to that for the 2D, we need to remove the errors $x_m$.  This can be done somewhat trivially compared with the sinogram alignment problem for parallel beam data, by using simple conservation of mass ideas.  That is to say, $\int Rf(x,y,\theta_m) \, dy \approx M(x)$, where $M(x)$ is the total mass along the $x$ cross section independent of $\theta_m$.  This approach has been successfully used in several instances, and can also be used to correct for offset or unknown tilt axis rotations as well \cite{wang2017fast,sanders2015physically,sanders2017mm}. If this condition is not sufficient for some reason, it may alternatively be achieved by manually alignment based on simple inspection.  Once this error has been removed, then for each cross-section $x$ we essentially have the 2D alignment problem considered throughout this paper.  Now however, the $\epsilon$ errors are the same for each cross-section.  Therefore we have hundreds or even thousands of 2D cross sections to drive the PBA, and the algorithm can be adapted in a natural way.

\subsection{Relationship to Projection Matching and New Improvement for Projection Matching}\label{sec:pm}

\begin{figure}
 \centering
 \includegraphics[width=.5\textwidth]{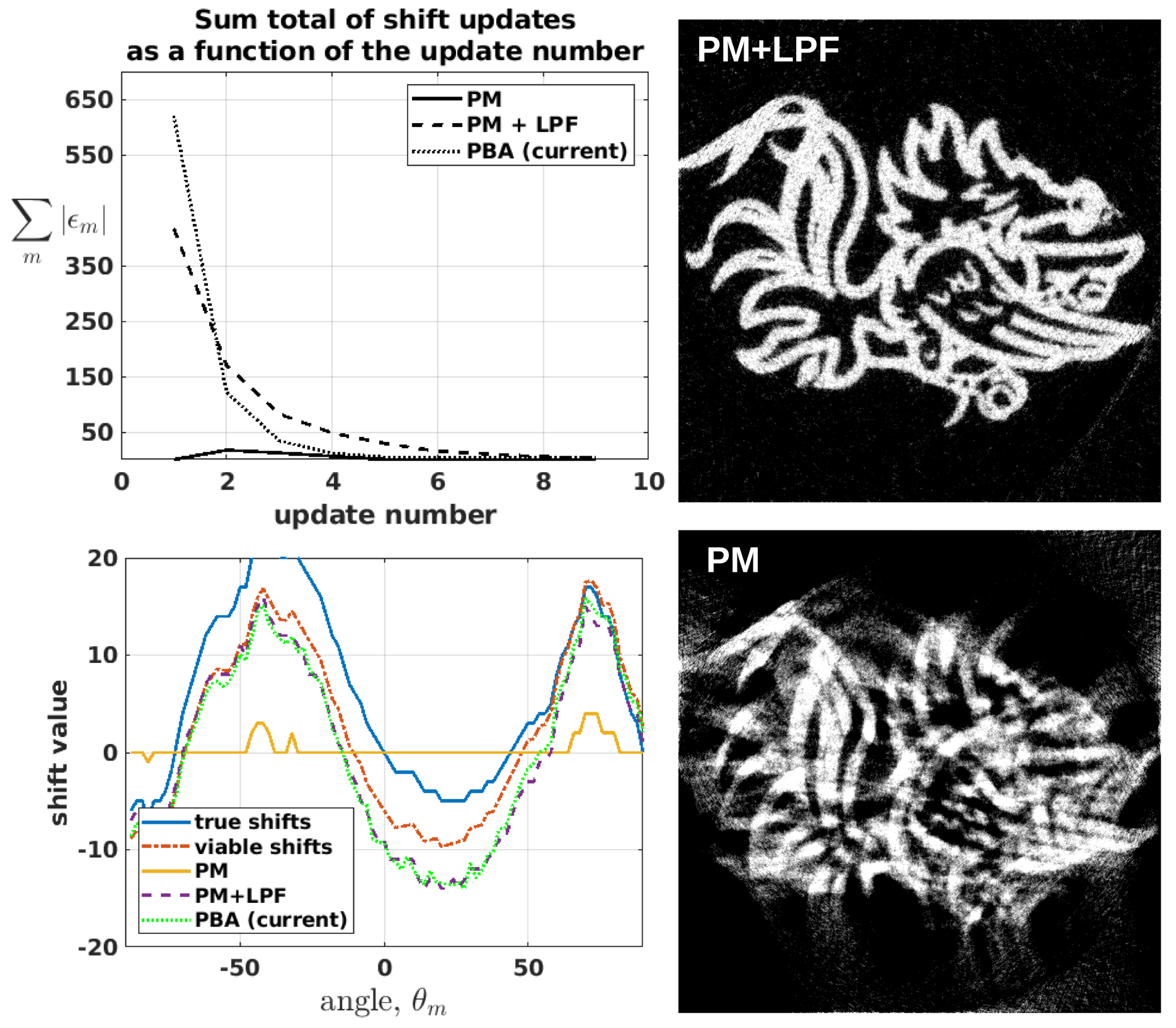}
 \caption{Repeating the cross correlation set up from Figure \ref{fig3}, we compare our method (PBA) with the PM and the PM with our proposed low pass filtering (LPF).  We observe that PM with LPF gives results very close to our method.}
 \label{fig:pm}
\end{figure}

\begin{table}
\begin{center}
\begin{tabular}{l | c | c | c | c|}
\hline
       Gamecock phantom errors       &  PBA & PM+LPF & PM & None \\
       \hline
       CC misaligned & .2999 & .3295 & .7327 & .7441 \\
       random misaligned (mean 0) & .2950  & .3216 & .5541 & .7431 \\
       random misaligned (mean 20) & .3998 & .4167 & .7081 & .8080 \\
       \hline
       Shepp-Logan errors & & & & \\
       \hline 
       CC misaligned & .2938 & .3836 & .7096 & .7096 \\
       random misaligned (mean 0) & .2948 & .3307 & .6450 & .7764 \\
       random misaligned (mean 20) & .2982 & .3237 & .7821 & .8446 \\
       \hline
\end{tabular}
\caption{Reconstruction errors from two test images dependent on the source of misalignment and the alignment method.  The PBA and PM+LPF are very similar in 5 out of 6 cases, and the PM alone is significantly worse in all cases.}
\end{center}
\label{table1}
\end{table}

Some time after the initial development of our alignment technique, we were able to develop a clear relationship to the iterative alignment method which we refer to as projection matching (PM) \cite{dengler1989multi,penczek1994ribosome,parkinson2012automatic}.  First both approaches work by iteratively refining the alignment based on discrepancies between the current solution and the data.  The gist of PM is that the acquired misaligned projections are iteratively realigned by registering/matching, with the corresponding projection with the current solution.  More formally, let the misaligned projection again be denoted by $\widetilde{Rf}(x,\theta_m)$, and let the projection at angle $\theta_m$ of the current solution $f^0$ be given by $Rf^0(x,\theta_m)$.  In its simplest form, the updated alignment (characterized in the continous setting) is given by cross correlating $\widetilde{Rf}(x,\theta_m)$ with $Rf^0(x,\theta_m)$, i.e. to find the new alignment shift $\epsilon_m$ by
\begin{equation}\label{proj-match}
 \max_{\epsilon_m} \int \widetilde{Rf}(x + \epsilon_m,\theta_m) Rf^0(x,\theta_m) \, dx.
\end{equation}
This approach has been shown to provide modest improvements, however the refined alignment tends to show slow convergence \cite{houben2011refinement}.  In our experience, any refinement is only realized with this approach whenever the data is severely misaligned, and therefore the refinement does not converge to an accurate final alignment.

To connect this approach more closely with our method, taking Fourier transforms and using Parseval's theorem, we see that (\ref{proj-match}) is equivalent to
\begin{equation}\label{pm-fourier}
 \max_{\epsilon_m} \int e^{i k \epsilon_m}  F(\widetilde{Rf} (\cdot , \theta_m ) )_k F( Rf^0(\cdot,\theta_m))_k \, dk.
\end{equation}
Here $F$ is taking the continuous 1D Fourier transform along the $x$ variable, analogous to our previous description of taking the DFT's, $F(r_{\theta_m})$.  With this reformulation the approach looks more like our method.  However, our method includes only the low frequency values that uses a precise formula to determine $\epsilon_m$.  Therefore, PM would relate more closely to our approach if a low pass filter (LPF) was applied before the matching, or alternatively only computing the integrand in (\ref{pm-fourier}) over small values $k$.  Figure \ref{fig:pm} demonstrates this idea, where we have repeated the set up from Figure \ref{fig3} with cross correlation as the misalignment and add noise with a SNR of 15, and applied PM with and without the LPF.  As shown in the recovered shift values in the bottom left panel, the PM without LPF yields essentially no or very slow shifts corrections, and with LPF the method recovers shifts values that are nearly identical to our phase based approach.  The convergence of the shifts are shown in the top left panel, where PM with LPF converges only slightly slower than our approach but could likely be improved with fine parameter tuning.  
This connection provides insight into why PM by cross-correlation is more robust when combined with a LPF.
In addition, PM is more seamlessly extended to the 3D setting.

The same simulation was repeated for the Shepp-Logan phantom with $5\degree$ angle increments, and the reconstruction errors were measured for all cases, including both random misalignment and misalignment due to cross correlation.  For the random misalignment, we include the cases of random misalignment shifts chosen on $[-20,20]$, and a nonzero mean set of misalignment shifts on the interval $[0,40]$.  The second set of shifts is reported as mean 20, and has the effect of offsetting the axis.  These errors are reported in Table \ref{table1}, and show that our approach closely matches the results of PM with an LPF, while PM by cross-correlation alone in its simplest form yields inferior results that are only a slight improvement to no alignment.

We briefly mention here that some algorithms with the flavor of PM have extended further than what mentioned here, by including additional alignment parameters and also formulating the problem with the maximum likelihood interpretation \cite{brandt2007structure,scheres2012bayesian}.  The alignment parameters include correction for an unstable tilt axis and inaccurate angle measures.  Recently in \cite{gursoy2017rapid}, some improvements were made to increase the convergence of projection matching, however as we observe here the PM combined with LPF provides significantly increased convergence in the recovery of the shift values.

\subsection{Non-uniqueness of the Correct Data Registration}\label{sec:nonu}
At first glance, it may seem that there is some unique solution for the optimal shifts.  This however is not true, and theoretically infinitely many solutions can exist.  This can be explained in simple terms by the fact that two unique sets of recovered shifts may recover two solutions which are precisely the same \emph{after} a translation of the recovered images, i.e. we may recover $\tilde f(x,y) = f(x-\alpha , y - \beta)$.  The main purpose of pointing this out is for the error analysis used in this paper.  For one, we obviously cannot simply measure the recovered shifts verses some applied misalignment shifts.  We also cannot simply measure the error in the reconstruction from the true test image.  We instead first determine the translation to match the reconstructed image with the test image, and then measure the relative error. In this domain it is trivial to see that the translation should be determined by cross-correlation, since cross-correlation of two images yields the shift which minimizes the distance (or relative error) between the images.

To formally show how this translation of the image results in the sinogram, we define  $\xi_\theta = (\cos \theta ,\sin \theta)$ and $t_{\alpha , \beta }^\theta  = t - (\alpha , \beta ) \cdot \xi_\theta $ and make use of the equivalent definition of the Radon transform, $Rf(t,\theta) = \iint_{\R^2} f(x,y) \delta(t-(x,y)\cdot \xi_\theta) \, dy \, dx$.  With these definitions, using a change of variables we can show $R \tilde f (t,\theta ) = Rf(t_{\alpha , \beta}^\theta,\theta).$
Essentially this informs us that the translated image by $(\alpha , \beta)$ results in shifts of the columns of the sinogram by $(\alpha , \beta)\cdot \xi_\theta$. 


\section{Further Simulations}
In this section we present results from various other simulations to demonstrate the robustness of our approach. We implement our method with the addition of TV regularization and 4 additional test images.  For the TV algorithm, the PBA is naturally meshed with the alternating direction method of multipliers (ADMM) algorithm, which is openly available at \cite{toby-web}.  Included in the supplementary material, we provide additional simulations demonstrating our method's robustness to limited data.

\subsection{Implementation with TV Regularization and 5 Test Images}

\begin{figure*}
\centering
 \includegraphics[width=1\textwidth]{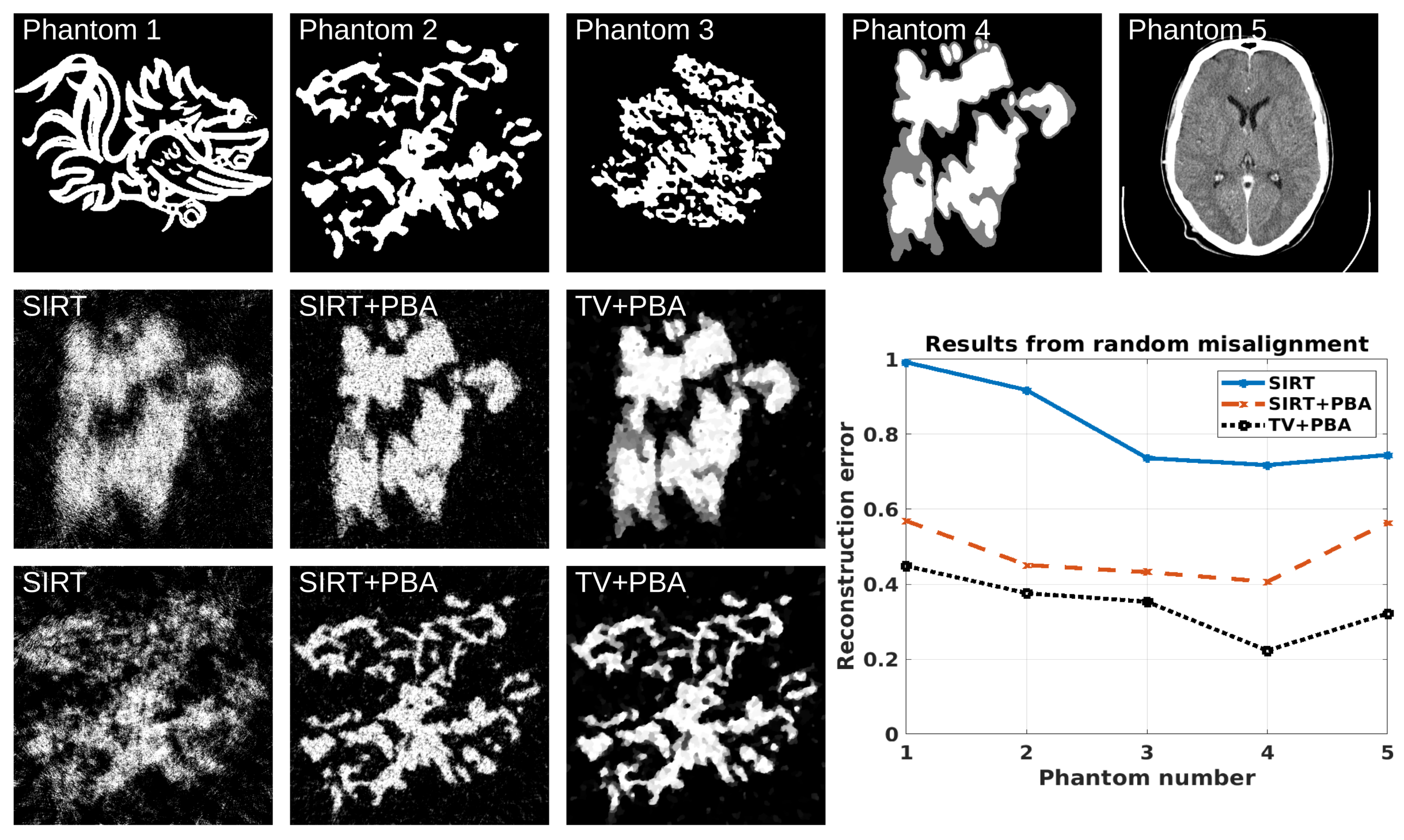}
 \caption{PBA results for 5 phantoms.  Test phantom images are shown in the top tow.  Reconstructions for phantom 4 and 5 are shown in the middle and bottom rows respectively.  The reconstruction errors are plotted in the bottom right panel.}
 \label{fig4}
\end{figure*}

Figure \ref{fig4} shows the 5 test images in the top row that we use for our simulations.  For these simulations, we simply apply random uniformly distributed misalignment shifts and apply our PBA algorithm as before.  The data again was generated over the full $180\degree$ range at every $2\degree$ and have added independent and identically distrubuted (i.i.d.) mean zero Gaussian noise with the variance set so that the signal to noise ratio (SNR) is 5.    The reconstruction results for phantoms 3 and 5 are shown in the bottom left 2 rows.  These results are fairly straight forward and appear to yield near exact alignment with the PBA both when using TV and SIRT.  The final error between the reconstructions and the original phantoms (after a correlation, necessary as pointed out in \ref{sec:nonu}) are given in the plot in the bottom right.  Obviously the PBA reconstructions yield significantly better errors than without, and TV also indicates some modest improvement here.  Generally speaking, these results indicate that our approach is effective with both SIRT and TV on a wide variety of images.

\subsection{Noise Analysis}
\begin{figure}
 \centering
 \includegraphics[width=.5\textwidth]{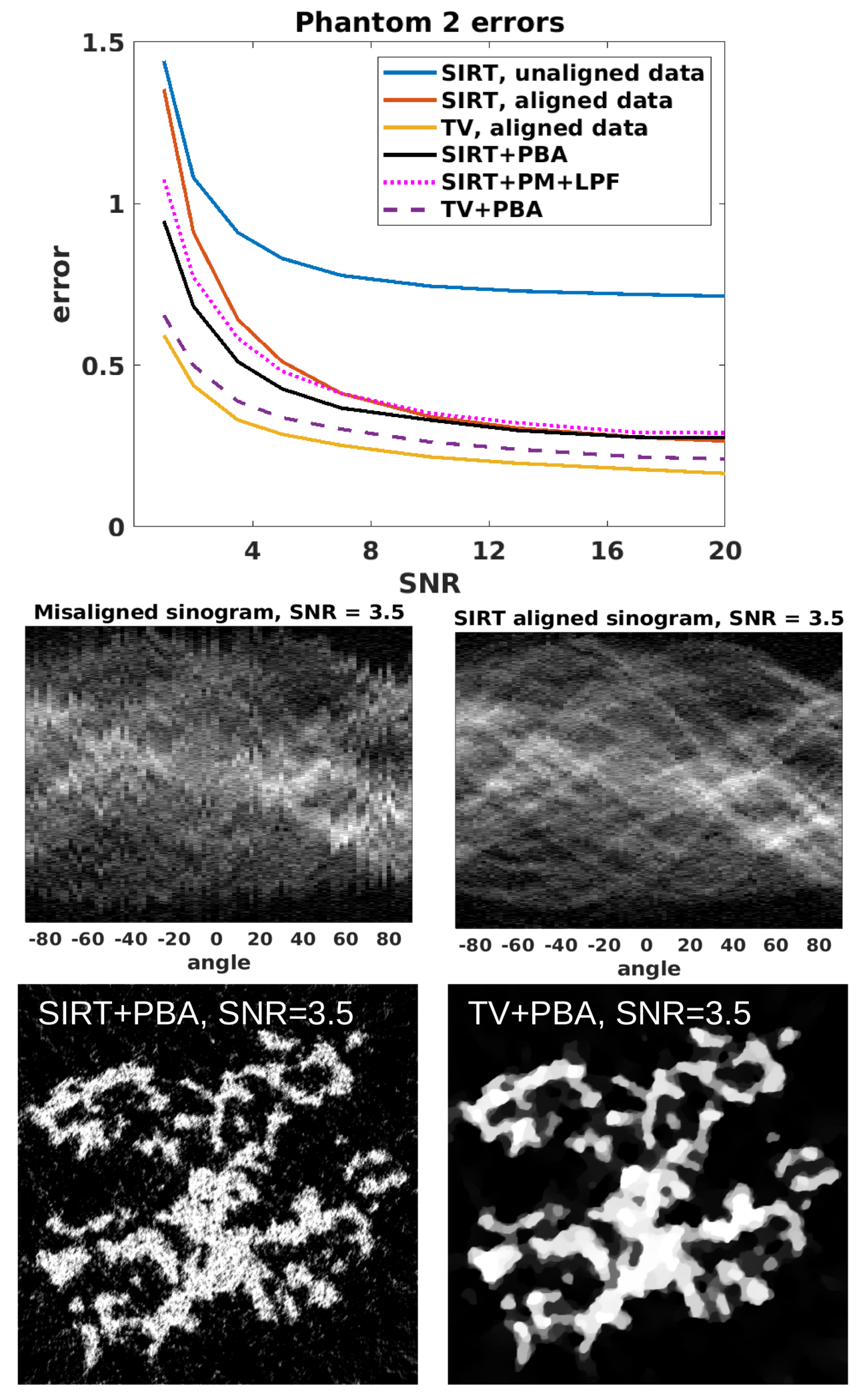}
 \caption{Reconstruction errors of the various alignment methods as a function of the SNR.  The misaligned and realigned sinogram are shown for phantom 2 at an SNR of 3.5, along with the PBA reconstructions from SIRT and TV.}
 \label{fig: noise}
\end{figure}

Here we present analysis on the robustness of our method with respect to increasingly high levels of noise.  For each of the 5 test images in Figure \ref{fig4}, we applied random misalignment and i.i.d. mean zero Gaussian noise to the sinograms, where the variance was set to yield a desired SNR.  We then applied the various reconstruction and alignment algorithms to this data.  The resulting reconstruction errors for phantom 2 are shown in Figure \ref{fig: noise} along with a few sample sinograms and reconstruction images\footnote{The general results were very similar for the other 4 phantom images and are thus provided in the supplementary work}.  For comparison purposes, we also included the reconstructions from noisy but perfectly aligned data to show the best results we could hope for in the alignment.  We see clearly here that TV provides a notable improvement especially at low SNR values.  We also again observe that PM with LPF yields results very similar to our proposed method, and in all cases these alignment methods yield approximately the same error as with the perfectly aligned data (red and yellow curves)\footnote{The SIRT reconstructions are actually slightly improved compared with perfectly aligned data, an artifact of overall iteration count with noisy data and ill-conditioned problems \cite{andersen1984simultaneous}.}.  Hence our method is achieving nearly the best possible results at each SNR level.

\section{Results on Electron Tomography Data Set}

\begin{figure}
\centering
 \includegraphics[width=.5\textwidth]{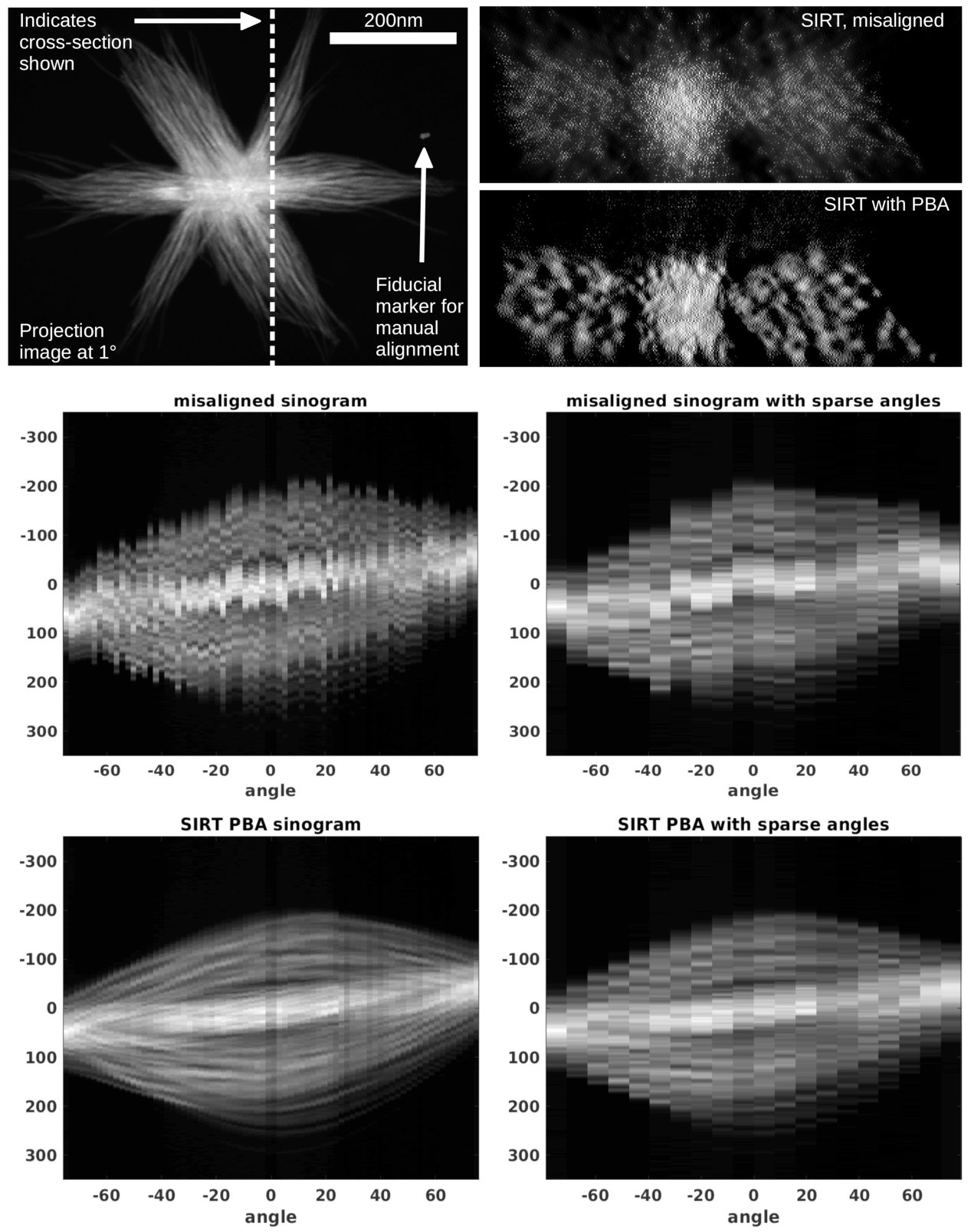}
 \caption{Alignment of electron microscopy tomographic tilt series.}
 \label{exp1}
\end{figure}

\begin{figure*}
 \centering
 \includegraphics[width=1\textwidth]{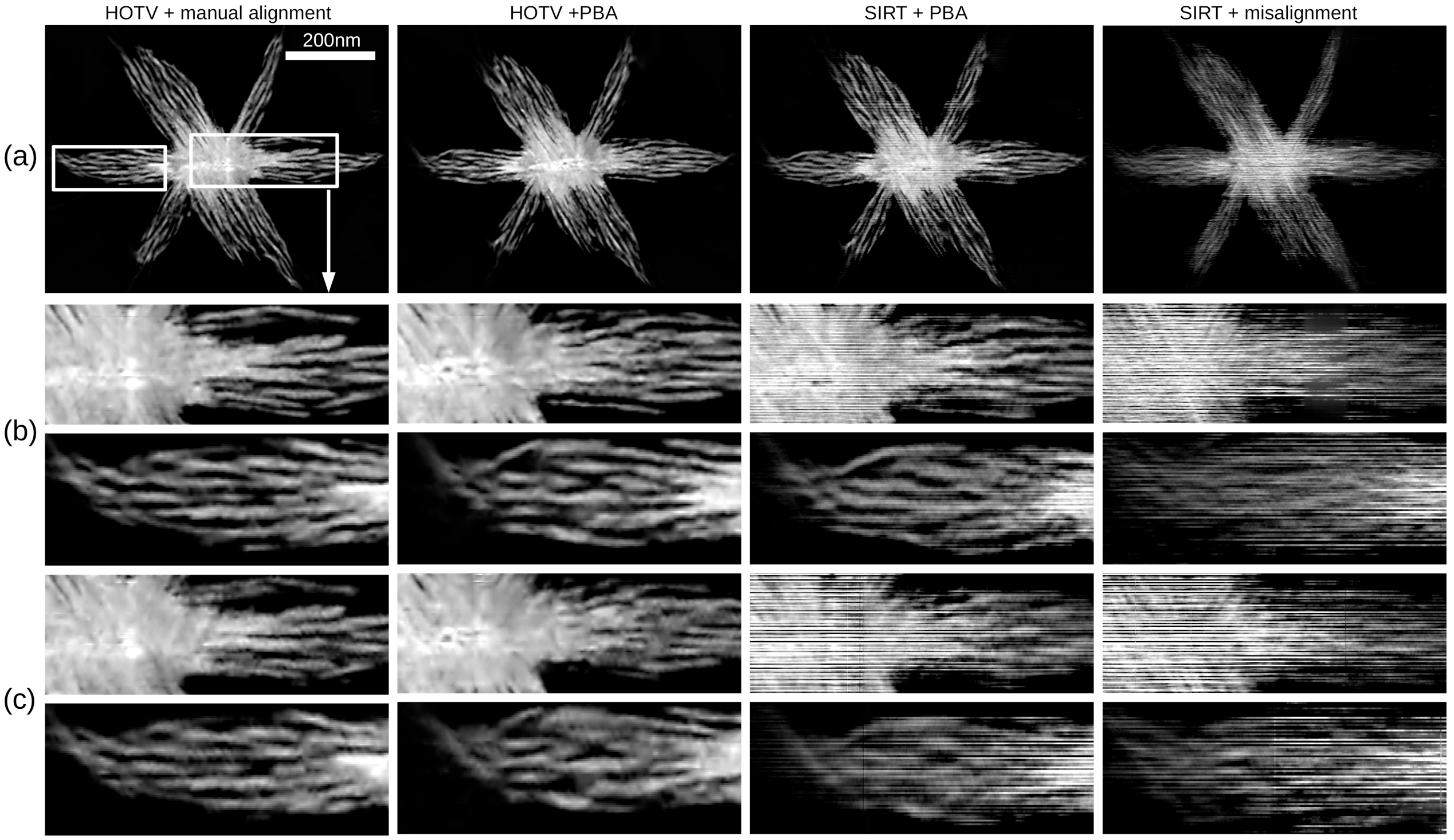}
 \caption{Cross-sectional images of the 3D electron micrscopy reconstructions from the various alignments and reconstruction algorithms.  In (a) are the results from the full data set ($2\degree$ angle increments) and (b) shows magnified patches from (a).  In (c) are the same magnified patches after reducing the data to $8\degree$ angle increments.}
 \label{exp2}
\end{figure*}

Electron tomography is a type of tomography where 3D nanoscale images are reconstructed from a series of 2D projection images.  The projection images are acquired with electron microscopes, and due to the small scale of the rotating mechanisms, alignment of the data is always necessary.  The ideal aligned version of the data are given by (\ref{3D-radon}), hence similar to the 2D version except there's and extra dimension, which is the axis of rotation.

Our data is provided via a recently published open access data project \cite{levin2016nanomaterial}.  The projections are acquired by rotating around a single access, and the relative angles are given at $\pm 75 \degree$ at every $2\degree$ (although there are two missing angles that can be observed at $23\degree$ and $25\degree$).  The authors of this data set aligned the projections manually based on the location of a single fiducial marker. Obviously, one would like to avoid the manual alignment if possible.  For our experiments, we introduced random misalignments back in the data, and naturally extended our method to 3D as described in Section \ref{sec:3D}, which we describe in finer detail here.  To challenge our methodology futher, we ran a second set of experiments by removing $75\%$ of the data and increasing the angle increment from $2\degree$ to $8\degree$.

Compared with the 2D simulations, here we have a very similar problem, however with this data set we have 1157 possible cross-sections in order to drive the PBA.  An accurate 3D reconstruction of this size may require several hours.  Therefore for practical implementation of these ideas, it makes more sense to implement the PBA over some subsets of cross-sections first, so that we do not potentially waste hours of computation for a poor result.  Therefore 20 cross-sections were selected first, and the phase based alignment was applied based on our general algorithm.  The shifts were determined by averaging the result of the estimated shifts for each cross-section, and allowing for 4 total iterates of alignment updates.  These recovered shifts were determined to provide a suitable alignment, and then the full 3D reconstruction was implemented.  If these shifts were determined unsuitable for any reason, the full set of 1157 could be utilized as an alternative.  The same approach was taken in the limited data experiment as well.

The various aligned sinograms from the full data (left) and the limited data (right) of a single cross-section are shown in Figure \ref{exp1}, along with a single projection image of the object and reconstructions of the cross-section using SIRT.  We have indicated the cross-section shown within the projection image, as well as the fiducial marker that the authors used for the alignment.  In the sinograms, it is difficult to make any solid conclusions on the quality of the alignment without running a reconstruction.  The two SIRT cross-sections shown here show significant improvement with the PBA and recover very sharp features.

More careful examination of the reconstruction is shown in Figure \ref{exp2}, where the cross-section shown is directly through the bulk of the particle, orthogonal to the cross-section displayed in Figure \ref{exp1}.  For the $\ell_1$ regularized solutions, we used higher order TV of order two, which is known to be more accurate than TV for these types of structures \cite{sanders2016}. Small patches of the images are magnified in the lower panels (b) and (c), which are indicated in the top-left image.  The images in (c) are from the limited data experiments, and (a) and (b) are from the full data.  All of the results clearly indicate an accurate alignment from the PBA compared with the misaligned data.  From a general inspection, it is difficult to determine whether the originally manually aligned data is more accurate than the PBA.  Indeed, all of these reconstruction images indicate that the data are well aligned, however not all of the structures appear precisely the same across the various methods.  One positive indicator that we can observe in the top row of panel (b) are the pores that are visible in the PBA image directly in the center of the object, which are not visible with the manually aligned images.  These pores can also vaguely been seen in the limited data case with HOTV.  Most other differences in the solutions are somewhat subjective.  Nonetheless, our automated PBA can achieve as good or possibly even better results than what was done with the manual alignment, and obviously significanly improves from the result of the unaligned data.

\section{Summary}
We have presented a new phase based method for alignment of sinogram data for accurate tomographic reconstructions.  The PBA works so that the data alignment is iteratively determined and automatically implemented within the numerical reconstruction algorithm.  The approach works by converting the misalignment physical shifts into multiplicative phase shifts in the data's Fourier domain, and the phase shifts are accurately estimated by the discrepancies between the misaligned data and a current reconstructed solution.  We have demonstrated that this approach is extremely accurate and robust in the variety of cases considered.  Namely, it is robust to limited data and noisy data, and it is effective with both SIRT and TV algorithms.  We also showed that our algorithm performed well with an electron microscopy data set, and even showed some modest improvement in the alignment compared with the manually aligned data provided from the original authors \cite{levin2016nanomaterial}.  While our approach shows great promise, future studies should be carried out to realize the full potential of our method.  These may also include adding additional misalignment parameters, such as inaccurate angle measurements.  Nonetheless, we feel confident in the ability of this approach for the advancement of automated tomographic imaging.


\end{document}